\documentstyle[11pt]{article}
\begin{document}
\newcommand{\qed}{\hphantom{.}\hfill $\Box$\medbreak}
\newcommand{\proof}{\noindent{\bf Proof \ }}
\newtheorem{Theorem}{Theorem}[section]
\newtheorem{Lemma}[Theorem]{Lemma}
\newtheorem{Corollary}[Theorem]{Corollary}
\newtheorem{Remark}[Theorem]{Remark}
\newtheorem{Example}[Theorem]{Example}
\newtheorem{Definition}[Theorem]{Definition}
\newtheorem{Construction}[Theorem]{Construction}

\thispagestyle{empty}

\begin{center}
{\Large\bf Enumerations of $(K_4-e)$-designs with small orders
\footnote{This work was supported in part by NSFC grant No.
$61071221$ (Y. Chang), NSFC grant No. $10901016$ (T. Feng) and by
PRIN, PRA and INDAM (GNSAGA) (G. Lo Faro and A. Tripodi).}}

\vskip12pt

Yanxun Chang, Tao Feng\\ Institute of Mathematics\\ Beijing Jiaotong
University\\ Beijing 100044, P. R. China\\ {\tt
yxchang@bjtu.edu.cn}\\ {\tt tfeng@bjtu.edu.cn}

\vskip6pt

Giovanni Lo Faro and Antoinette Tripodi\\ Department of
Mathematics\\ University of Messina\\ Viale Ferdinando Stagno
d'Alcontres, 31 - 98166, Sant'Agata\\ Messina, Italy\\ {\tt lofaro@unime.it}\\
{\tt tripodi@dipmat.unime.it}
\end{center}

\begin{center}
{\it Dedicated to the memory of Lucia Gionfriddo $(1973$-$2008)$}
\end{center}

\vskip8pt

\noindent {\bf Abstract:} It is established that up to isomorphism,
there are only one $(K_4-e)$-design of order $6$, three
$(K_4-e)$-designs of order $10$ and two $(K_4-e)$-designs of order
$11$. As an application of our enumerative results, we discuss the
fine triangle intersection problem for $(K_4-e)$-designs of orders
$v=6,10,11$.



\section{Introduction}
Let $H$ be a simple graph and $G$ be a subgraph of $H$. A {\em
$G$-design} of $H$ (or $(H, G)$-design) is a pair $(X, {\cal B})$
where $X$ is the vertex set of $H$ and ${\cal B}$ is an
edge-disjoint decomposition of $H$ into isomorphic copies (called
{\em blocks}) of the graph $G$. If $H$ is the complete graph $K_v$,
we refer to such a $G$-design as one of order $v$.

The most basic question in design theory is that given a graph $G$
and a positive integer $v$, whether a $G$-design of order $v$
exists. If the existence problem is answered completely, then a
further question is what about the enumeration problem for
$G$-designs of order $v$. That is to say up to isomorphism, how many
$G$-designs of order $v$ exist? Two $G$-designs of order $v$
$(X,{\cal B}_1)$ and $(X,{\cal B}_2)$ are said to be {\em
isomorphic} if there exists a permutation $\pi$ on $X$ such that
$\pi({\cal B}_1)={\cal B}_2$, where $\pi$ is applied to the elements
of each block of ${\cal B}_1$. For more information on $G$-designs,
the interested reader may refer to \cite{abb}.

If $G$ is the graph with vertices $a$, $b$, $c$, $d$ and edges $ab$,
$ac$, $ad$, $bc$, $bd$, then such a graph is called a $K_4-e$ and
denoted by $[a,b,c-d]$. Bermond and Sch\"{o}nheim established that
\cite{bs} a $(K_4-e)$-design of order $v$ exists if and only if
$v\equiv 0,1\ ({\rm mod }\ 5)$ and $v\geq 6$. In this paper we will
focus on the enumerations of $(K_4-e)$-designs of order $v$ for
$v=6,10,11$. We shall show that there is only one $(K_4-e)$-design
of order $6$ up to isomorphism. There are three non-isomorphic
$(K_4-e)$-designs of order $10$ and two non-isomorphic
$(K_4-e)$-designs of order $11$.

Finally, as an application of our enumerative results, we
investigate the fine triangle intersection problem for
$(K_4-e)$-designs of orders $v=6,10,11$.

\section{Enumerations of $v=6$ and $10$}

\begin{Theorem}
\label{6non-isomorphic} There is only one $(K_4-e)$-design of order $6$ up to isomorphism.
\end{Theorem}

\proof Let $X=\{0,1,2,3,4,5\}$. Suppose that $(X,{\cal B})$ is any
$(K_4-e)$-design of order $6$. For every $x\in X$, denote by
$d_i(x)$, $i=2,3$, the number of blocks of ${\cal B}$ in which the
degree of $x$ is $i$. It follows that $2d_2(x)+3d_3(x)=5$, and so
$d_2(x)=d_3(x)=1$. It is readily checked that up to isomorphism the
unique $(K_4-e)$-design of order $6$ is:

\begin{center}
\begin{tabular}{llll}
\hspace{3cm}${\cal B}:$ & $[0,1,2-3]$,& $[2,3,4-5]$,& $[4,5,0-1]$.
\end{tabular}\qed
\end{center}

\begin{Lemma}
\label{10lemma1} Any $(K_4-e)$-design of order $10$ contains a
subdesign of order $6$.
\end{Lemma}

\proof Let $(X,{\cal B})$ be any $(K_4-e)$-design of order $10$. For
every $x \in X$, denote by $d_i(x)$, $i=2,3$, the number of blocks
of ${\cal B}$ in which the degree of $x$ is $i$. It follows that
$2d_2(x)+3d_3(x)=9$. Solving this equation gives two possibilities:
$d_2(x)=0$ and $d_3(x)=3$ (we refer to such a vertex as a
$a$-element) or $d_2(x)=3$ and $d_3(x)=1$ (we refer to such a vertex
as a $b$-element). Denote the number of $a$-elements and
$b$-elements by $\alpha$ and $\beta$, respectively. Since each block
contains exactly two elements with degree 3 we have
$$\left\{ \begin{array}{l}
                  3\alpha +\beta =18\\
                  \alpha +\beta =10
                  \end{array}
          \right. $$
and so $\alpha =4$ and $\beta =6$. Let $A=\{6,7,8,9\}$ and
$B=\{0,1,2,3,4,5\}$ be the sets of $a$-elements and $b$-elements,
respectively. Then we have the following six blocks ${\cal
C}=\{[6,7,x_1-x_2]$, $[8,9,x_1-x_2]$, $[6,8,x_3-x_4]$,
$[7,9,x_3-x_4]$, $[6,9,x_5-x_6]$, $[7,8,x_5-x_6]\}$, where $\{x_1,$
$x_2,$ $x_3,$ $x_4,$ $x_5,$ $x_6\}=\{0,$ $1, 2, 3, 4, 5\}$, and so
$(B,{\cal B}\setminus {\cal C})$ is a $(K_4-e)$-design of order
6.\qed

\begin{Theorem}
\label{10non-isomorphic} There are exactly $3$ non-isomorphic
$(K_4-e)$-designs of order $10$.
\end{Theorem}

\proof Let $X=\{0, 1, 2, \ldots, 9\}$. Suppose that $(X,{\cal B})$
is any $(K_4-e)$-design of order 10. Note that there is the only
unique $(K_4-e)$-design of order $6$ under isomorphism. By Lemma
\ref{10lemma1}, $(X,{\cal B})$ contains a subdesign of order $6$,
say, $[0,1,2-3]$, $[2,3,4-5]$, $[4,5,0-1]$. The other $6$ blocks of
${\cal B}$ must be the following forms: $[6,7,x_1-x_2]$,
$[8,9,x_1-x_2]$, $[6,8,x_3-x_4]$, $[7,9,x_3-x_4]$, $[6,9,x_5-x_6]$,
$[7,8,x_5-x_6]\}$, where $\{x_1,$ $x_2,$ $x_3,$ $x_4,$ $x_5,$
$x_6\}=\{0, 1, 2, 3, 4, 5\}$.

Let $D=\{ \{0,1\}, \{2,3\}, \{4,5\} \}$ and $N=\{ \{x_1, x_2\},
\{x_3, x_4\}, \{x_5, x_6\} \}$. We distinguish the following three
cases.

(1) If $|D\cap N|=3$, without loss of generality we can always
assume that $\{x_1, x_2\}=\{0,1\}$. We then have $\{x_3,
x_4\}=\{2,3\}$ and $\{x_5, x_6\}=\{4,5\}$ under isomorphism.

(2) If $|D\cap N|=1$, similarly we can always assume that $\{x_1,
x_2\}=\{0,1\}$. We then have $\{x_3, x_4\}=\{2,4\}$ and $\{x_5,
x_6\}=\{3,5\}$ under isomorphism.

(3) If $|D\cap N|=0$, similarly we can always assume that $\{x_1,
x_2\}=\{0,2\}$. We then have $\{x_3, x_4\}=\{1,4\}$ and $\{x_5,
x_6\}=\{3,5\}$ under isomorphism.

From the above discussions, we have the following three
$(K_4-e)$-designs of order $10$ under isomorphism:
\begin{center}
  \begin{tabular}{llllll}
${\cal B}_1:$ & $[0,1,2-3]$, & $[2,3,4-5]$, & $[4,5,0-1]$, &
$[6,7,0-1]$, & $[8,9,0-1]$,\\
& $[6,8,2-3]$, & $[7,9,2-3]$, & $[6,9,4-5]$, & $[7,8,4-5]$;\\
\\

${\cal B}_2:$ & $[0,1,2-3]$, & $[2,3,4-5]$, & $[4,5,0-1]$, &
$[6,7,0-1]$, & $[8,9,0-1]$,\\
& $[6,8,2-4]$, & $[7,9,2-4]$, & $[6,9,3-5]$, & $[7,8,3-5]$;\\
\\

${\cal B}_3:$ & $[0,1,2-3]$, & $[2,3,4-5]$, & $[4,5,0-1]$, &
$[6,7,0-2]$, & $[8,9,0-2]$,\\
& $[6,8,1-4]$, & $[7,9,1-4]$, & $[6,9,3-5]$, & $[7,8,3-5]$.
\end{tabular}
\end{center}
It is easy to see that any two of the three $(K_4-e)$-designs of
order $10$ are not isomorphic. This completes the proof. \qed

\section{Enumeration of $v=11$}

\begin{Lemma}
\label{v=11 D-cycle} Let $(X, {\cal B})$ be any $(K_4-e)$-design of
order $11$. Let $D=\{ \{x,y\}:\ [x,y,z-u]\in {\cal B}\}$. Then $(X,
D)$ is a cycle of length $11$.
\end{Lemma}

\proof For every $x \in X$, denote by $d_i(x)$, $i=2,3$, the number
of blocks of ${\cal B}$ in which the degree of $x$ is $i$. It
follows that $2d_2(x)+3d_3(x)=10$, which gives two solutions:
$d_2(x)=2$ and $d_3(x)=2$ (we refer to such a vertex as a
$a$-element) or $d_2(x)=5$ and $d_3(x)=0$ (we refer to such a vertex
as a $b$-element). Denote the number of $a$-elements and
$b$-elements by $\alpha$ and $\beta$, respectively. Since each block
contains exactly two elements with degree 3 we have $2\alpha =22$
and $\alpha +\beta =11$, and so $\alpha =11$ and $\beta =0$. We then
have the fact that for every $x \in X$
\begin{equation}\label{fact-11}
d_2(x)=d_3(x)=2.
\end{equation}
By the definition of $D$, $(X,D)$ is a $2$-regular graph, whose
connected component is a cycle with length at least $3$. Let
$X=\{0,1,2,\ldots,10\}$. We will show that $(X,D)$ is a cycle of
length $11$. All possibilities are exhausted as below.

(1). If $(X,D)$ contains a cycle of length $4$, without loss of
generality, we can assume that $\cal B$ contain the following four
blocks:  $[0,1,x_1,x_2]$, $[1,2,x_3,x_4]$, $[2,3,x_5,x_6]$,
$[3,0,x_7,x_8]$, where $x_1, x_2, \ldots, x_8\in X$. From the fact
(\ref{fact-11}), $\{0,2\}$ must appear in some one of the four
above-listed blocks. Since $(X,{\cal B})$ is a $(K_4-e)$-design of
order $11$, we have $x_i\not\in \{0,1,2,3\}$ for each $1\leq i\leq
8$. It is impossible.

(2). If $(X,D)$ contains a cycle of length $5$, similarly we can
assume that $\cal B$ contain the following five blocks:
$[0,1,x_1,x_2]$, $[1,2,x_3,x_4]$, $[2,3,x_5,x_6]$, $[3,4,x_7,x_8]$,
$[4,0,x_9,x_{10}]$, where $x_1, x_2, \ldots, x_{10}\in X$. From the
fact (\ref{fact-11}), the $5$ remaining $2$-subsets of
$\{0,1,\ldots,4\}$, $\{0,2\}$, $\{0,3\}$, $\{1,3\}$, $\{1,4\}$,
$\{2,4\}$, must appear in some one of the five above-listed blocks.
Note that $|\{x_{2i-1}, x_{2i}\}\cap \{0,1,\ldots,4\}|\leq 1$ for
each $i=1,2,\ldots,5$. If $|\{x_{2i-1}, x_{2i}\}\cap
\{0,1,\ldots,4\}|=1$ for some $i$, then there appear exactly two
$2$-subsets of $\{ \{0,2\}, \{0,3\}, \{1,3\}, \{1,4\}, \{2,4\} \}$
in the block containing $\{x_{2i-1}, x_{2i}\}$. A contradiction
occurs.

(3). If $(X,D)$ contains a cycle of length $3$, similarly we can
assume that $\cal B$ contain the following three blocks:
$[0,1,3-4]$, $[1,2,5-6]$, $[2,0,7-8]$. Since $\{0,x\}$, $\{1,x\}$
and $\{2,x\}$ where $x\in\{9,10\}$ must appear in some block of
$\cal B$, by the fact (\ref{fact-11}) we can further assume that
$\cal B$ contain the following blocks: $[9,10,0-1]$, $[9,c,2-*]$,
$[10,d,2-*]$, $[e,f, 0-*]$, $[g,h,1-*]$ under isomorphism. Since the
unused $2$-subsets containing $2$ are only $\{2,3\}$ and $\{2,4\}$.
So, we have $c=3$ and $d=4$, or $c=4$ and $d=3$. Under isomorphism,
we can assume that $c=3$ and $d=4$. Similar arguments give
$\{e,f\}=\{5,6\}$ and $\{g,h\}=\{7,8\}$. From the fact
(\ref{fact-11}), other three blocks of $\cal B$ can be assumed in
the two possibilities under isomorphism:

Case I: \ $[3,4, *-*]$, $[5,7,*-*]$, $[6,8,*-*]$;

Case II: \ $[3,5, *-*]$, $[4,7,*-*]$, $[6,8,*-*]$.

For Case I, $(X,D)$ contains a cycle of length $4$: $\{9,10\}$,
$\{10,4\}$, $\{4,3\}$ and $\{3,9\}$. By the arguments of (2), it is
impossible.

For Case II, an exhausted search by hand shows that it is impossible
to complete the partial design to a $(K_4-e)$-design of order $11$.

This completes the proof. \qed

\begin{Lemma}
\label{v=11 N-cycle} Let $(X, {\cal B})$ be any $(K_4-e)$-design of
order $11$. Let $D=\{ \{x,y\}:\ [x,y,z-u]\in {\cal B}\}$ and $N=\{
\{z,u\}:\ [x,y,z-u]\in {\cal B}\}$. Then $D\cap N=\emptyset$.
\end{Lemma}

\proof Let $X=\{0,1,2,\ldots,10\}$. Assume that $D\cap N\neq\emptyset$. Without loss of generality let $\{\{0,1\}\}\in D\cap N$ and

\hspace{2cm} ${\cal B}_0=\{[0,1,2-3],[4,5,0-1],[0,6,7-8],[9,10,0-x_1]\}\subset {\cal B}$.

\noindent ${\cal B}_0$ is the set of all blocks containing $0$.

Note that by the formula $(1)$ in the proof of Lemma \ref{v=11
D-cycle}, for each $x\in X$ and each $i=2,3$, there are exactly $2$
blocks in $\cal B$ containing $x$ in which the degree of $x$ is $i$.
This fact will be used OFTEN in the following!

Furthermore, consider the blocks containing $1$. It is readily
checked that up to isomorphism, except for the blocks $[0,1,2-3]$
and $[4,5,0-1]$, the blocks containing $1$ must be one of the
following four cases:

\hspace{2cm} $(1)$ $[6,1,9-10]$, $[7,8,1-x_2]$;

\hspace{2cm} $(2)$ $[7,1,8-9]$, $[6,10,1-x_2]$;

\hspace{2cm} $(3)$ $[9,1,6-7]$, $[8,10,1-x_2]$;

\hspace{2cm} $(4)$ $[9,1,7-8]$, $[6,10,1-x_2]$.

\underline{Case $(1)$}. Let ${\cal B}_1=\{[6,1,9-10],[7,8,1-x_2]\}$.
Then ${\cal B}_0\cup{\cal B}_1\subset {\cal B}$. It follows that
$\{\{0,1\},\{1,6\},\{6,0\}\}\subset D$, which is a cycle of length
$3$. By Lemma \ref{v=11 D-cycle}, a contradiction.

\underline{Case $(2)$}. Let ${\cal B}_1=\{[7,1,8-9],[6,10,1-x_2]\}$.
Then ${\cal B}_0\cup{\cal B}_1\subset {\cal B}$. Consider the blocks
containing $6,7,8$, respectively. We have

\hspace{2cm}${\cal B}_2=\{[*,*,6-7],[*,*,6-*],[*,8,*-*],[*,8,*-*],[*,7,*-*]\}\subset {\cal B}$,

\noindent and ${\cal B}={\cal B}_0\cup{\cal B}_1\cup{\cal B}_2$.
Since the edge $\{7,10\}$ must occur in one block of $\cal B$, there
must be a block of the form $[*,7,10-*]$. Since the edge $\{8,10\}$
must occur in one block of $\cal B$, there must be a block of the
form $[*,8,10-*]$. Since the edge $\{6,9\}$ must occur in one block
of $\cal B$, there must be a block of the form $[*,9,6-*]$.
Combining the above three facts, we rewrite ${\cal B}_2$ as follows:

\hspace{2cm}${\cal B}_2=\{[*,*,6-7],[*,9,6-*],[*,8,10-*],[*,8,*-*],[*,7,10-*]\}$.

\noindent Consider the blocks containing the edge $\{8,9\}$. We have
the following two possibilities for ${\cal B}_2$:

Subcase I. \ $[x_3,x_4,6-7],[x_5,9,6-x_6],[x_7,8,10-9],[x_8,8,x_9,x_{10}],
[x_{11},7,10-x_{12}]$;

Subcase II. \ $[x_3,x_4,6-7],[x_5,9,6-x_6],[x_7,8,10-x_8],[x_9,8,9,x_{10}],
[x_{11},7,10-x_{12}]$.

\noindent It is readily checked that for each $1\leq i\leq 12$,
$x_i\in \{2,3,4,5\}$. Let $\{a,b,c,d\}=\{2,3,4,5\}$. Consider the
blocks containing the edges $\{2,10\},\{3,10\},\{4,10\},\{5,10\}$.
Then we can always take $x_1=a$, $x_2=b$, $x_7=c$, $x_{11}=d$.

For Subcase I, consider the blocks containing $9$. We have $x_5=d$
and $x_6=b$. Consider the blocks containing $6$. We have $x_3=a$ and
$x_4=c$. Consider the blocks containing $7$. We have $x_{12}=b$.
Thus $x_2=x_6=x_{12}=b$, and there are $3$ blocks in $\cal B$
containing $b$ in which the degree of $b$ is $2$. A contradiction.

For Subcase II, consider the blocks containing the edge $\{6,a\}$.
We have $x_3=a$. Consider the block containing $7$. We have $x_4=c$
and $x_{12}=b$. Consider the block containing $6$. We have $x_5=d$.
Consider the blocks containing $9$. We have $x_6=c$ and $x_9=b$.
Consider the blocks containing $8$. We have
$\{x_8,x_{10}\}=\{a,d\}$. If $x_8=a$ and $x_{10}=d$, then the edge
$\{a,c\}$ occurs in two blocks of $\cal B$. A contradiction. If
$x_8=d$ and $x_{10}=a$, then the edge $\{c,d\}$ occurs in two blocks
of $\cal B$. A contradiction.

\underline{Case $(3)$}. Let ${\cal B}_1=\{[9,1,6-7],[8,10,1-x_2]\}$.
Then ${\cal B}_0\cup{\cal B}_1\subset {\cal B}$. Consider the blocks
containing $6,7,9$, respectively. We have

\hspace{2cm}${\cal B}_2=\{[*,7,*-*],[*,7,*-*],[*,6,*-*],[*,*,6-9],[*,*,9-*]\}\subset {\cal B}$,

\noindent and ${\cal B}={\cal B}_0\cup{\cal B}_1\cup{\cal B}_2$.
Since the edge $\{6,10\}$ must occur in one block of $\cal B$, there
must be a block of the form $[*,6,10-*]$. Since the edge $\{7,10\}$
must occur in one block of $\cal B$, there must be a block of the
form $[*,7,10-*]$. Since the edge $\{8,9\}$ must occur in one block
of $\cal B$, there must be a block of the form $[*,8,9-*]$. Consider
the blocks containing $8$. Except the blocks $[0,6,7-8]$,
$[8,10,1-x_2]$ and $[*,8,9-*]$, there must be a block of the form
$[*,7,8-*]$. Combining the above four facts, we have the following
two possibilities for ${\cal B}_2$:

Subcase I. \ $[x_3,7,10-8],[x_4,7,x_5-x_6],[x_7,6,10-x_8],[x_9,x_{10},6-9],
[x_{11},8,9-x_{12}]$;

Subcase II. \ $[x_3,7,10-x_4],[x_5,7,8-x_6],[x_7,6,10-x_8],[x_9,x_{10},6-9],
[x_{11},8,9-x_{12}]$.

\noindent It is readily checked that for each $1\leq i\leq 12$,
$x_i\in \{2,3,4,5\}$. Let $\{a,b,c,d\}=\{2,3,4,5\}$. Consider the
blocks containing the edges $\{2,10\},\{3,10\},\{4,10\},\{5,10\}$.
Then we can always take $x_1=a$, $x_2=b$, $x_3=c$, $x_7=d$. Consider
the blocks containing $9$. We have $x_9=b$, $x_{10}=c$, $x_{11}=d$.
Consider the blocks containing $6$. We have $x_8=a$.

For Subcase I, consider the blocks containing $8$. We have
$x_{12}=a$.  Thus $x_1=x_8=x_{12}=a$, and there are $3$ blocks in
$\cal B$ containing $a$ in which the degree of $a$ is $2$. A
contradiction.

For Subcase II, consider the blocks containing $8$. We have $x_5=a$
and $x_{12}=c$. Consider the blocks containing $7$. We have
$\{x_4,x_6\}=\{b,d\}$. If $x_4=b$ and $x_6=d$, then the edge
$\{b,c\}$ occurs in two blocks. A contradiction. If $x_4=d$ and
$x_6=b$, then the edge $\{c,d\}$ occurs in two blocks. A
contradiction.

\underline{Case $(4)$}. Let ${\cal B}_1=\{[9,1,7-8],[6,10,1-x_2]\}$.
Then ${\cal B}_0\cup{\cal B}_1\subset {\cal B}$. It follows that
$\{\{0,1\},\{1,9\},\{9,10\},\{10,6\},\{6,0\}\}\subset D$, which is a
cycle of length $5$. By Lemma \ref{v=11 D-cycle}, a contradiction.

This completes the proof. \qed

\begin{Theorem}
\label{11non-isomorphic} There are exactly $2$ non-isomorphic $(K_4-e)$-designs of
order $11$.
\end{Theorem}

\proof Let $X=\{0,1,2,\ldots,10\}$. Suppose that $(X, {\cal B})$ is
any $(K_4-e)$-design of order $11$. Let $D=\{ \{x,y\}:\ [x,y,z-u]\in
{\cal B}\}$ and $N=\{ \{z,u\}:\ [x,y,z-u]\in {\cal B}\}$. By Lemma
\ref{v=11 D-cycle}, $(X, D)$ is a cycle of length $11$. Without loss
of generality, assume that $D=\{\{i,i+1\}:0\leq i\leq
9\}\cup\{\{10,0\}\}$. Consider the blocks $[10,0,a-b]$ and
$[0,1,c-d]$. We have $\{a,b\}\subset\{2,3,4,5,6,7,8\}$,
$\{c,d\}\subset\{3,4,5,6,7,8,9\}$, and
$\{a,b\}\cap\{c,d\}=\emptyset$. By Lemma \ref{v=11 N-cycle}, $D\cap
N=\emptyset$. Due to $\{\{a,b\},\{c,d\}\}\subset N$, it follows that
$\{a,b\}\in\{\{j,k\}:2\leq j\leq 6,j+2\leq k\leq 8\}$ and
$\{c,d\}\in\{\{j,k\}:3\leq j\leq 7,j+2\leq k\leq 9\}$.

By the formula $(1)$  in the proof of Lemma \ref{v=11 D-cycle}, for
every $x\in X$, there are exactly $2$ blocks in $\cal B$ containing
$x$ and satisfying the degree of $x$ is $2$. Let $[s_1,s_1+1,0-s_2]$
and $[t_1,t_1+1,0-t_2]$ be the two blocks containing $0$ in which
the degree of $0$ is $2$. Consider the edges containing $0$. We have
$\{10,1,a,b,c,d,s_1,s_1+1,t_1,t_1+1\}=X\setminus\{0\}$. An
exhaustive search by hand shows that there are $19$ possibilities
for the values $a,b,c,d,s_1,t_1$ satisfying

$(1)$ $\{a,b\}\in\{\{j,k\}:2\leq j\leq 6,j+2\leq k\leq 8\}$;

$(2)$ $\{c,d\}\in\{\{j,k\}:3\leq j\leq 7,j+2\leq k\leq 9\}$;

$(3)$ $\{10,1,a,b,c,d,s_1,s_1+1,t_1,t_1+1\}=X\setminus\{0\}$.

\noindent We list them in the first five columns in Table I.

\begin{center}
{\bf Table I}

\tabcolsep 0.08in
  \begin{tabular}{|c|c|c|c|c|c|c|}\hline
&$\{a,b\}$ & $\{c,d\}$ & $(s_1,s_1+1)$ & $(t_1,t_1+1)$ & $\pi\{a,b\}$ & $\pi\{c,d\}$\\\hline
$1^*$& $\{2,4\}$ &$\{3,5\}$ &$(6,7)$ &$(8,9)$ &$\{7,9\}$ &$\{6,8\}$\\
$2^*$& $\{2,4\}$ &$\{3,7\}$ &$(5,6)$ &$(8,9)$ &$\{7,9\}$ &$\{4,8\}$\\
$3^*$& $\{2,4\}$ &$\{3,9\}$ &$(5,6)$ &$(7,8)$ &$\{7,9\}$ &$\{2,8\}$\\\hline
$4^*$& $\{2,5\}$ &$\{6,9\}$ &$(3,4)$ &$(7,8)$ &$\{6,9\}$ &$\{2,5\}$\\\hline
$5^*$& $\{2,6\}$ &$\{3,7\}$ &$(4,5)$ &$(8,9)$ &$\{5,9\}$ &$\{4,8\}$\\
$6^*$& $\{2,6\}$ &$\{3,9\}$ &$(4,5)$ &$(7,8)$ &$\{5,9\}$ &$\{2,8\}$\\
$7^*$& $\{2,6\}$ &$\{5,7\}$ &$(3,4)$ &$(8,9)$ &$\{5,9\}$ &$\{4,6\}$\\
$8^*$& $\{2,6\}$ &$\{5,9\}$ &$(3,4)$ &$(7,8)$ &$\{5,9\}$ &$\{2,6\}$\\\hline
$9^*$& $\{2,7\}$ &$\{3,6\}$ &$(4,5)$ &$(8,9)$ &$\{4,9\}$ &$\{5,8\}$\\\hline
$10^*$& $\{2,8\}$ &$\{3,9\}$ &$(4,5)$ &$(6,7)$ &$\{3,9\}$ &$\{2,8\}$\\
$11$& $\{2,8\}$ &$\{7,9\}$ &$(3,4)$ &$(5,6)$ &$\{3,9\}$ &$\{2,4\}$\\
$12$& $\{2,8\}$ &$\{5,9\}$ &$(3,4)$ &$(6,7)$ &$\{3,9\}$ &$\{2,6\}$\\\hline
$13^*$& $\{4,6\}$ &$\{5,7\}$ &$(2,3)$ &$(8,9)$ &$\{5,7\}$ &$\{4,6\}$\\
$14$& $\{4,6\}$ &$\{5,9\}$ &$(2,3)$ &$(7,8)$ &$\{5,7\}$ &$\{2,6\}$\\\hline
$15$& $\{4,8\}$ &$\{5,9\}$ &$(2,3)$ &$(6,7)$ &$\{3,7\}$ &$\{2,6\}$\\
$16$& $\{4,8\}$ &$\{7,9\}$ &$(2,3)$ &$(5,6)$ &$\{3,7\}$ &$\{2,4\}$\\\hline
$17^*$& $\{5,7\}$ &$\{4,6\}$ &$(2,3)$ &$(8,9)$ &$\{4,6\}$ &$\{5,7\}$\\\hline
$18$& $\{5,8\}$ &$\{4,9\}$ &$(2,3)$ &$(6,7)$ &$\{3,6\}$ &$\{2,7\}$\\\hline
$19$& $\{6,8\}$ &$\{7,9\}$ &$(2,3)$ &$(4,5)$ &$\{3,5\}$ &$\{2,4\}$\\
\hline
\end{tabular}
\end{center}

Let $\pi=(0)(1\ 10)(2\ 9)(3\ 8)(4\ 7)(5\ 6)$ be a permutation on
$X$. Obviously $\pi(D)=D$. Since
$\pi([10,0,a-b])=[0,1,\pi(a)-\pi(b)]$ and
$\pi([0,1,c-d])=[10,0,\pi(c)-\pi(d)]$, under the action of $\pi$, if
for some possibility in Table I, whose values in the second and
third column are $\{a,b\}$ and $\{c,d\}$, respectively, then it is
isomorphic to the possibility with values $\pi\{c,d\}$ and
$\pi\{a,b\}$ in the second and third column, respectively. Using
this idea the above $19$ possibilities can be reduced to $12$
possibilities. We mark them with a $*$. Take the first possibility
for example. In the first possibility, $\{a,b\}=\{2,4\}$ and
$\{c,d\}=\{3,5\}$. Then $\pi\{c,d\}=\{6,8\}$ and
$\pi\{a,b\}=\{7,9\}$, which corresponds to the last possibility.

Next consider the blocks $[0,1,c-d]$ and $[1,2,e-f]$. For each given
$\{c,d\}$ in Table I, we need to determine all possible values of
$\{e,f\}$. Fix the permutation $\sigma=(0\ 1\ 2\ 3\ 4\ 5\ 6\ 7\ 8\
9\ 10)$ on $X$ and let $A$ be the collection of the set
$\{[10,0,a-b],[0,1,c-d]\}$, where $\{a,b\}$, $\{c,d\}$ are taken
from the $19$ possibilities in Table I. Let $B$ be the collection of
all possible cases of the set $\{[0,1,c_1-d_1],[1,2,e-f]\}$ such
that one can complete a $(K_4-e)$-design originally from the two
blocks $[0,1,c_1-d_1]$ and $[1,2,e-f]\}$. It is easy to see that
$B\subseteq \sigma(A)$. Thus for determining $B$, we count
$\sigma(A)$. Apply permutation $\sigma$ to Table I to obtain Table
II. 
Note that because of
$\{10,1,a,b,c,d,s_1,s_1+1,t_1,t_1+1\}=X\setminus\{0\}$, we have
$\sigma(\{10,1,a,b,c,d,s_1,s_1+1,t_1,t_1+1\})=X\setminus\{1\}$. Here 
$[s'_1,s'_1+1,1-s'_2]$ and $[t'_1,t'_1+1,1-t'_2]$ are the two blocks
containing $1$ in which the degree of $1$ is $2$.

\begin{center}
{\bf Table II}

\tabcolsep 0.08in
  \begin{tabular}{|c|c|c|c|c|}\hline
&$\{c_1,d_1\}$ & $\{e,f\}$ & $(s'_1,s'_1+1)$ & $(t'_1,t'_1+1)$\\\hline
$1$& $\{3,5\}$ &$\{4,6\}$ &$(7,8)$ &$(9,10)$\\
$2$& $\{3,5\}$ &$\{4,8\}$ &$(6,7)$ &$(9,10)$\\
$3$& $\{3,5\}$ &$\{4,10\}$ &$(6,7)$ &$(8,9)$\\\hline
$4$& $\{3,6\}$ &$\{7,10\}$ &$(4,5)$ &$(8,9)$\\\hline
$5$& $\{3,7\}$ &$\{4,8\}$ &$(5,6)$ &$(9,10)$\\
$6$& $\{3,7\}$ &$\{4,10\}$ &$(5,6)$ &$(8,9)$\\
$7$& $\{3,7\}$ &$\{6,8\}$ &$(4,5)$ &$(9,10)$\\
$8$& $\{3,7\}$ &$\{6,10\}$ &$(4,5)$ &$(8,9)$\\\hline
$9$& $\{3,8\}$ &$\{4,7\}$ &$(5,6)$ &$(9,10)$\\\hline
$10$& $\{3,9\}$ &$\{4,10\}$ &$(5,6)$ &$(7,8)$\\
$11$& $\{3,9\}$ &$\{8,10\}$ &$(4,5)$ &$(6,7)$\\
$12$& $\{3,9\}$ &$\{6,10\}$ &$(4,5)$ &$(7,8)$\\\hline
$13$& $\{5,7\}$ &$\{6,8\}$ &$(3,4)$ &$(9,10)$\\
$14$& $\{5,7\}$ &$\{6,10\}$ &$(3,4)$ &$(8,9)$\\\hline
$15$& $\{5,9\}$ &$\{6,10\}$ &$(3,4)$ &$(7,8)$\\
$16$& $\{5,9\}$ &$\{8,10\}$ &$(3,4)$ &$(6,7)$\\\hline
$17$& $\{6,8\}$ &$\{5,7\}$ &$(3,4)$ &$(9,10)$\\\hline
$18$& $\{6,9\}$ &$\{5,10\}$ &$(3,4)$ &$(7,8)$\\\hline
$19$& $\{7,9\}$ &$\{8,10\}$ &$(3,4)$ &$(5,6)$\\\hline
\end{tabular}
\end{center}

Now for each given $\{a,b\}$ and $\{c,d\}$ in Table I, we can use
Table II to determine all possible values of $\{e,f\}$. For example
when we take the first possibility in Table I, i.e.,
$\{a,b\}=\{2,4\}$ and $\{c,d\}=\{3,5\}$, the values of $(e,f)$ can
be taken from the rows with $(c_1,d_1)=(3,5)$ in Table II. Thus we
have $(e,f)=(4,6)$, $(4,8)$ or $(4,10)$. These three subcases
corresponding to the first possibility in Table I are listed below.

\begin{center}
{\bf Table III}

\tabcolsep 0.08in
  \begin{tabular}{|c|c|c|c|c|c|c|c|}\hline
&$\{a,b\}$ & $\{c,d\}$ & $\{e,f\}$ & $(s_1,s_1+1)$ & $(t_1,t_1+1)$ & $(s'_1,s'_1+1)$ & $(t'_1,t'_1+1)$\\\hline
$1^*$& $\{2,4\}$ &$\{3,5\}$ &$\{4,6\}$ & $(6,7)$ &$(8,9)$  &$(7,8)$ &$(9,10)$\\
& & & $\{4,8\}$ &&&$(6,7)$ &$(9,10)$\\
& & & $\{4,10\}$ &&&$(6,7)$ &$(8,9)$\\\hline
\end{tabular}
\end{center}

For reducing these subcases, we notice that $s_1$ can not be equal
to $s'_1$. Otherwise there would be a block $[s_1,s_1+1,0-1]\in{\cal
B}$, which implies $\{0,1\}\in N$. Due to $\{0,1\}\in D$ and $D\cap
N=\emptyset$ from Lemma \ref{v=11 N-cycle}, a contradiction occurs.
Similarly, we have $s_1\neq t'_1$, $t_1\neq s'_1$ and $t_1\neq
t'_1$. Thus $|\{s_1,t_1,s'_1,t'_1\}|=4$. Using this condition, for
each given $\{a,b\}$ and $\{c,d\}$ in Table I, we can reduce
possible values of $\{e,f\}$. For example in Table III only the
first subcase satisfies $|\{s_1,t_1,s'_1,t'_1\}|=4$. After
exhaustive search by hand, we can reduce the $12$ possibilities
marked $^*$ in Table I to $6$ possibilities in Table IV.

\begin{center}
{\bf Table IV}

\tabcolsep 0.08in
  \begin{tabular}{|c|c|c|c|c|c|c|c|}\hline
&$\{a,b\}$ & $\{c,d\}$ & $\{e,f\}$ & $(s_1,s_1+1)$ & $(t_1,t_1+1)$ & $(s'_1,s'_1+1)$ & $(t'_1,t'_1+1)$\\\hline
$1^*$& $\{2,4\}$ &$\{3,5\}$ &$\{4,6\}$ & $(6,7)$ &$(8,9)$  &$(7,8)$ &$(9,10)$\\\hline
$2^*$& $\{2,4\}$ &$\{3,7\}$ &$\{6,8\}$ & $(5,6)$ &$(8,9)$  &$(4,5)$ &$(9,10)$\\\hline
$3^*$& $\{2,4\}$ &$\{3,9\}$ &$\{8,10\}$ & $(5,6)$ &$(7,8)$  &$(4,5)$ &$(6,7)$\\\hline
$5^*$& $\{2,6\}$ &$\{3,7\}$ &$\{4,8\}$ & $(4,5)$ &$(8,9)$  &$(5,6)$ &$(9,10)$\\\hline
$10^*$& $\{2,8\}$ &$\{3,9\}$ &$\{4,10\}$ & $(4,5)$ &$(6,7)$  &$(5,6)$ &$(7,8)$\\\hline
$13^*$& $\{4,6\}$ &$\{5,7\}$ &$\{6,8\}$ & $(2,3)$ &$(8,9)$  &$(3,4)$ &$(9,10)$\\\hline
\end{tabular}
\end{center}

In the following we show that Possibilities $2$, $3$, $10$ and $13$
in Table IV are impossible. For Possibilities $2$ and $3$, consider
the block containing the edge $\{2,4\}$. It must be of the form
$[2,3,4-*]$ (note that $N\cap D=\emptyset$). Since there are four
blocks containing $2$ in $\cal B$, and it is easy to verify that one
can not find the fourth block containing $2$, a contradiction
occurs. For Possibility $10$, the edge $\{2,10\}$ occurs in two
blocks $[10,0,2-8]$ and $[1,2,4-10]$. A contradiction. For
Possibility $13$, consider the block containing the edge $\{4,6\}$.
It must be of the form $[6,7,4-*]$. Consider the block containing
the edge $\{6,8\}$. It must be of the form $[5,6,8-*]$. Then all
blocks containing $4$ are of the form $[10,0,4-6]$, $[6,7,4-*]$,
$[3,4,*-*]$, $[4,5,*-*]$, and all blocks containing $8$ are of the
form $[1,2,6-8]$, $[5,6,8-*]$, $[7,8,*-*]$, $[8,9,*-*]$. It follows
that there is no block containing the edge $\{4,8\}$. A
contradiction.

By Possibilities $1$ and $5$ in Table IV, the blocks $[10,0,a-b]$
and $[0,1,c-d]$ must be one of the following two cases: $(1)$
$[10,0,2-4]$, $[0,1,3-5]$; $(2)$ $[10,0,2-6]$, $[0,1,3-7]$. It
implies that for any $i\in X$, the blocks $[i,i+1,x-y]$ and
$[i+1,i+2,z-u]$ must be one of the following two cases: $(1)$
$[i,i+1,(i+3)-(i+5)]$, $[i+1,i+2,(i+4)-(i+6)]$; $(2)$
$[i,i+1,(i+3)-(i+7)]$, $[i+1,i+2,(i+4)-(i+8)]$, where the arithmetic
is modulo $11$. It follows that $\cal B$ must be one of the
following two cases:
\begin{center}
${\cal B}_1=\{[i,i+1,(i+3)-(i+5)]:i\in X\}$,

${\cal B}_2=\{[i,i+1,(i+3)-(i+7)]:i\in X\}$.
\end{center}

\noindent It is readily checked that $(X,{\cal B}_1)$ and $(X,{\cal
B}_2)$ are both $(K_4-e)$-designs, and they are non-isomorphic. This
completes the proof. \qed

\section{Application in fine triangle intersection problem}

As an application of our enumerative results, in this section we
investigate the fine triangle intersection problem for
$(K_4-e)$-designs of orders $v=6,10,11$.

Let $B$ be a simple graph. Denote by $T(B)$ the set of all triangles
of the graph $B$. For example, if $B=[a,b,c-d]$, then
$T(B)=\{\{a,b,c\},\{a,b,d\}\}$. Two $G$-designs of order $v$
$(X,{\cal B}_1)$ and $(X,{\cal B}_2)$ {\em intersect} in $t$
triangles provided $|T({\cal B}_1)\cap T({\cal B}_2)|=t$, where
$T({\cal B}_i)=\bigcup_{B\in {\cal B}_i} T(B)$, $i=1,2$. Define
$Fin_G(v)=\{(s,t):$ $\exists$ a pair of $G$-designs of order $v$
intersecting in $s$ blocks and $t+s|T(G)|$ triangles$\}$. The {\em
fine triangle intersection problem} for $G$-designs is to determine
$Fin_G(v)$.

The fine triangle intersection problem for $G$-designs, which was
introduced in \cite{cflt}, is the generalization of the intersection
problem and the triangle intersection problem for $G$-designs. For
more information on the intersection problem for $G$-designs, the
interested reader may refer to \cite{bgl,bk,cl,chl,km,lr}. For more
information on the triangle intersection problem for $G$-designs,
the interested reader may refer to \cite{bly,cfl,ly}.

Let $b_v=v(v-1)/10$ be the number of blocks in a $( K_4-e)$-design,
and $[a,b]$ the set of all integers $x$ satisfying $a\leq x \leq b$.
Let $J(v)=\{s:$ there exists a pair of $(K_4-e)$-designs of order
$v$ intersecting in $s$ blocks$\}$, and $J_T(v)=\{t:$ there exists a
pair of $(K_4-e)$-designs of order $v$ intersecting in $t$
triangles$\}$.

\begin{Theorem}
\label{intersection-(K_4-e)} {\rm (\cite{bgl})}  For any $v\equiv
0,1\ ({\rm mod }\ 5)$, $v\geq 6$ and $v\neq 11$, $J(v)=[0,
b_v]\setminus\{b_v-1, b_v-2\}$; $J(11)=\{0, 1, 2, \ldots, 6, 11\}.$
\end{Theorem}

\begin{Theorem}
\label{triangel-intersection for K_4-e} {\rm (\cite{bly})} For any
$v\equiv 0,1\ ({\rm mod }\ 5)$, $v\geq 15$,
$J_T(v)=[0,2b_v]\setminus\{2b_v-1,2b_v-2\}$; $J_T(6)=\{0,2,3,6\};$
$J_T(10)=\{0,1, \ldots, 12,14,15,18\};$ $J_T(11)=\{0,1, \ldots,
16,22\}$.
\end{Theorem}

If a pair of $(K_4-e)$-designs have blocks in common, each common
block contributes $2$ common triangles. In what follows we always
write $Fin_G(v)$ simply as $Fin(v)$ when $G$ is the graph $K_4-e$,
i.e., $Fin(v)=\{(s,t):$ $\exists$ a pair of $(K_4-e)$-designs of
order $v$ intersecting in $s$ blocks and $t+2s$ triangles$\}$. Let
$Adm(v)=\{(s,t):\ s+t\leq b_v, s\in J(v), \ t+2s\in J_T(v)\}\}$.
From the definitions of $Fin(v)$, $J(v)$ and $J_T(v)$, it is clear
that $Fin(v)\subseteq Adm(v)$.

\begin{Theorem}
\label{6} $Fin(6)=Adm(6)$.
\end{Theorem}
\proof Let $X=\{0,1,2,3,4,5\}$ and ${\cal B}=\{[0,1,2-3],[2,3,4-5],
[4,5,0-1]\}$. Then $(X,{\cal B})$ is a $(K_4-e)$-design of order
$6$. Consider the following permutations on $X$.

\begin{center}
  \begin{tabular}{llll}
$\pi_{0,0}=(2\ 4)(3\ 5)$, & $\pi_{0,2}=(1\ 2)$, & $\pi_{0,3}=(1\ 3)(2\ 4)$, & $\pi_{3,0}=(1)$.
\end{tabular}
\end{center}

\noindent Then we have $|\pi_{s,t}{\cal B}\cap {\cal
B}|=s$ and $|T(\pi_{s,t}{\cal B}\setminus{\cal B})\cap T({\cal
B}\setminus\pi_{s,t}{\cal B})|=t$ for each $(s,t)\in Adm(6)$. \qed

\begin{Theorem}
\label{10} $Fin(10)=Adm(10)\setminus\{(1,8),(3,1),(3,5),(4,1),(4,3),(5,1),(5,2)\}$.
\end{Theorem}
\proof Take the same $(K_4-e)$-designs of order $10$ $(X,{\cal
B}_i)$, $i=1,2,3$, as those in the proof of Theorem
\ref{10non-isomorphic}, which are mutually non-isomorphic. Consider
the following permutations on $X$.

\begin{center}\scriptsize
  \begin{tabular}{llll}
$\pi_{0,0}=(2\ 4)(3\ 5)(7\ 8)$,&
$\pi_{0,1}=(3\ 4\ 6)(5\ 8\ 9)$,&
$\pi_{0,2}=(3\ 4\ 5\ 6)(7\ 8)$,&
$\pi_{0,3}=(3\ 4)(5\ 6)(7\ 8)$,\\
$\pi_{0,4}=(3\ 4)(7\ 8)$,&
$\pi_{0,5}=(3\ 6)(5\ 8)$,&
$\pi_{0,6}=(1\ 2)(4\ 6)(5\ 7)$,&
$\pi_{0,7}=(1\ 2)(3\ 4)$,\\
$\pi_{0,8}=(1\ 2)(3\ 6)(5\ 7)$,&
$\pi_{0,9}=(1\ 2)(3\ 4)(8\ 9)$,&
$\pi_{1,0}=(4\ 6)(5\ 8\ 9)$,&
$\pi_{1,1}=(5\ 6)(7\ 8\ 9)$,\\
$\pi_{1,2}=(2\ 4)(3\ 6)$,&
$\pi_{1,3}=(5\ 6)(7\ 8)$,&
$\pi_{1,4}=(3\ 4)(5\ 6)$,&
$\pi_{1,5}=(3\ 6)(5\ 7)$,\\
$\pi_{1,6}=(1\ 6)(2\ 4\ 5\ 9)(3\ 7)$,&
$\pi_{1,7}=(1\ 6)(5\ 8)$,&
$\pi_{2,0}=(4\ 6)(5\ 7\ 8\ 9)$,&
$\pi_{2,1}=(5\ 6)(8\ 9)$,\\
$\pi_{2,2}=(4\ 6)(5\ 7\ 8)$,&
$\pi_{2,3}=(3\ 4\ 6)(5\ 7)$,&
$\pi_{2,4}=(3\ 6)(4\ 8)(5\ 9)$,&
$\pi_{2,5}=(3\ 5\ 7\ 4\ 6)(8\ 9)$,\\
$\pi_{2,6}=(3\ 4)$,&
$\pi_{2,7}=(1\ 2)(3\ 5)$,&
$\pi_{3,0}=(7\ 8\ 9)$,&
$\pi_{3,2}=(4\ 5)(7\ 8)$,\\
$\pi_{3,3}=(5\ 6)$,&
$\pi_{3,4}=(0\ 2\ 4)(1\ 3\ 5)$,&
$\pi_{3,6}=(0\ 2\ 4)(1\ 3\ 5)(7\ 8)$,&
$\pi_{4,0}=(4\ 6)(5\ 7)(8\ 9)$,\\
$\pi_{4,2}=(4\ 6)(5\ 7)$,&
$\pi_{4,4}=(4\ 6)(5\ 7)$,&
$\pi_{5,0}=(8\ 9)$,&
$\pi_{5,4}=(4\ 5)$,\\
$\pi_{6,0}=(2\ 4)(3\ 5)(8\ 9)$,&
$\pi_{6,2}=(3\ 5)$,&
$\pi_{6,3}=(0\ 2)(1\ 4)(7\ 8)$,&
$\pi_{9,0}=(1)$.
\end{tabular}
\end{center}

\noindent Let $E=\{(1,8),(3,1),(3,5),(4,1),(4,3),(5,1),(5,2)\}$ and $M=\{(1,6)$, $(1,7)$, $(2,5)$, $(2,7)$, $(3,2)$, $(3,4)$, $(3,6)$, $(4,0)$, $(4,4)$, $(5,4)$, $(6,2)$, $(6,3)\}$. Then for each $(s,t)\in Adm(10)\setminus(E\cup M)$, $|\pi_{s,t}{\cal B}_1\cap {\cal B}_1|=s$ and $|T(\pi_{s,t}{\cal B}_1\setminus{\cal B}_1)\cap T({\cal
B}_1\setminus\pi_{s,t}{\cal B}_1)|=t$. For each $(s,t)\in M$, $|\pi_{s,t}{\cal B}_2\cap {\cal B}_2|=s$ and $|T(\pi_{s,t}{\cal B}_2\setminus{\cal B}_2)\cap T({\cal
B}_2\setminus\pi_{s,t}{\cal B}_2)|=t$.

Now it remains to show that for each $(s,t)\in E$, we have $(s,t)\not\in Fin(10)$.
By Theorem \ref{10non-isomorphic}, there are exactly $3$ non-isomorphic $(K_4-e)$-designs of order $10$. Thus we can check all the cases by computer exhaustive search for the fine triangle intersection numbers of a pair of $(K_4-e)$-designs of order $10$, i.e., for any permutation $\pi$ on $X$ and for each $i,j=1,2,3$, count $|\pi{\cal B}_i\cap {\cal B}_j|$ and $|T(\pi{\cal B}_i\setminus{\cal B}_j)\cap T({\cal
B}_j\setminus\pi{\cal B}_i)|$. This completes the proof. \qed

\begin{Theorem}
\label{11} $Fin(11)=Adm(11)\setminus\{(3,0)$, $(4,0)$, $(4,6)$, $(5,0)$, $(5,1)$, $(5,2)$, $(6,0)$,  $(6,1)$, $(6,3)$, $(6,4)\}$.
\end{Theorem}
\proof Take $X=\{0,1,2,\ldots,10\}$ and ${\cal B}_1=\{[i,i+1,(i+3)-(i+5)]:i\in X\}$,
${\cal B}_2=\{[i,i+1,(i+3)-(i+7)]:i\in X\}$, where the arithmetic is modulo $11$. By Theorem \ref{11non-isomorphic}, $(X,{\cal B}_1)$ and $(X,{\cal B}_2)$ are two non-isomorphic $(K_4-e)$-designs of order $11$. Consider the following permutations on $X$.
\begin{center}\scriptsize
  \begin{tabular}{llll}
$\pi_{0,0}=(0\ 8\ 5\ 4\ 6)(2\ 9\ 10\ 3)$,&
$\pi_{0,1}=(0\ 8\ 1\ 7)(2\ 3\ 9\ 6\ 5\ 4)$,&
$\pi_{0,2}=(0\ 10\ 8\ 4)(2\ 9\ 6\ 7\ 3)$,\\
$\pi_{0,3}=(0\ 3)(2\ 8\ 7\ 4\ 5\ 10\ 6\ 9)$,&
$\pi_{0,4}=(2\ 6\ 8\ 3\ 5\ 9\ 10\ 7\ 4)$,&
$\pi_{0,5}=(0\ 9\ 4)(1\ 6\ 10\ 3)(2\ 5\ 8\ 7)$,\\
$\pi_{0,6}=(0\ 5\ 10\ 9)(1\ 7\ 2)(3\ 8)$,&
$\pi_{0,7}=(0\ 2\ 3\ 9\ 6)(1\ 8)(4\ 10)$,&
$\pi_{0,8}=(0\ 3\ 2\ 6\ 8)(1\ 5\ 4\ 9\ 7)$,\\
$\pi_{0,9}=(0\ 3\ 7\ 6\ 8\ 2)(1\ 10)(5\ 9)$,&
$\pi_{0,10}=(0\ 7\ 10\ 6\ 8\ 5\ 2\ 1\ 3\ 9\ 4)$,&
$\pi_{0,11}=(0\ 2\ 6\ 1\ 7\ 4)(5\ 9\ 8)$,\\
$\pi_{1,0}=(0\ 8\ 3\ 4\ 10)(1\ 2\ 5\ 7\ 9\ 6)$,&
$\pi_{1,1}=(0\ 2\ 7\ 3\ 8)(1\ 6\ 4\ 10\ 5\ 9)$,&
$\pi_{1,2}=(0\ 10\ 2\ 1\ 5)(3\ 7\ 4\ 6\ 9)$,\\
$\pi_{1,3}=(0\ 2\ 6\ 7\ 8\ 9)(1\ 4\ 3\ 5\ 10)$,&
$\pi_{1,4}=(0\ 1\ 7\ 5\ 10\ 9\ 6\ 3)(4\ 8)$,&
$\pi_{1,5}=(0\ 5\ 8\ 7)(1\ 2\ 4\ 6)(3\ 10\ 9)$,\\
$\pi_{1,6}=(0\ 3\ 7\ 1\ 8\ 10\ 5\ 4\ 2\ 9\ 6)$,&
$\pi_{1,7}=(0\ 5\ 3\ 10)(1\ 2\ 7\ 8\ 9\ 4)$,&
$\pi_{1,8}=(0\ 8\ 10\ 7\ 5\ 3)(1\ 6\ 4\ 2)$,\\
$\pi_{1,9}=(0\ 10\ 2\ 8\ 9\ 3\ 7\ 4)(1\ 5)$,&
$\pi_{1,10}=(1\ 10\ 3\ 4\ 9\ 6\ 5\ 2)(7\ 8)$,&
$\pi_{2,0}=(0\ 10\ 5\ 6)(1\ 7)(2\ 3\ 9\ 4\ 8)$,\\
$\pi_{2,1}=(0\ 8\ 1\ 7\ 6\ 3\ 9\ 5\ 2)(4\ 10)$,&
$\pi_{2,2}=(0\ 2\ 3\ 6\ 10\ 1)(4\ 5)(7\ 9)$,&
$\pi_{2,3}=(0\ 1\ 5\ 9\ 7\ 3\ 10\ 2\ 6\ 8)$,\\
$\pi_{2,4}=(1\ 8\ 5)(2\ 10\ 4)(6\ 7)$,&
$\pi_{2,5}=(0\ 10\ 4\ 5\ 9\ 3\ 8\ 1\ 6\ 7\ 2)$,&
$\pi_{2,6}=(0\ 7\ 3)(1\ 6\ 10\ 8\ 4\ 2)(5\ 9)$,\\
$\pi_{2,7}=(1\ 3\ 6\ 7\ 4\ 5)$,&
$\pi_{2,8}=(0\ 2\ 4\ 8)(1\ 3\ 7\ 6)(9\ 10)$,&
$\pi_{2,9}=(0\ 4\ 6)(1\ 9\ 2\ 8)(3\ 7\ 10)$,\\
$\pi_{3,1}=(0\ 9\ 7\ 5\ 3)(1\ 2)(4\ 6)(8\ 10)$,&
$\pi_{3,2}=(0\ 5\ 1\ 9\ 3\ 7\ 10\ 4)(2\ 8\ 6)$,&
$\pi_{3,3}=(0\ 7\ 3\ 10\ 6\ 4)(1\ 8\ 5)(2\ 9)$,\\
$\pi_{3,4}=(0\ 9\ 5)(1\ 8\ 2\ 7\ 3\ 6\ 10\ 4)$,&
$\pi_{3,5}=(0\ 10\ 4\ 2)(1\ 9\ 7\ 5\ 3)(6\ 8)$,&
$\pi_{3,6}=(1\ 10\ 8\ 6\ 9\ 7\ 5\ 2)(3\ 4)$,\\
$\pi_{3,7}=(0\ 2)(3\ 6)$,&
$\pi_{3,8}=(0\ 3\ 8\ 2\ 7\ 1\ 6)(4\ 9\ 5)$,&
   $\pi_{4,1}=(0\ 2\ 6\ 5\ 10\ 3\ 1\ 7\ 4)(8\ 9)$,\\
   $\pi_{4,2}=(0\ 8\ 9\ 4\ 5\ 6\ 7\ 10\ 3\ 1)$,&
$\pi_{4,3}=(0\ 5\ 9\ 1\ 4\ 8)(2\ 3\ 7\ 10)$,&
$\pi_{4,4}=(0\ 6\ 9\ 2\ 4\ 8\ 1\ 5\ 10\ 3\ 7)$,\\
$\pi_{4,5}=(0\ 3\ 6)(1\ 4\ 7\ 8\ 9)(2\ 5\ 10)$,&
  $\pi_{4,7}=(0\ 1\ 5\ 9\ 7\ 3)(2\ 6)(4\ 10)$,&
$\pi_{5,3}=(0\ 8\ 7\ 3\ 1\ 10\ 9\ 6\ 4\ 2)$,\\
$\pi_{5,4}=(0\ 10\ 9\ 8\ 7\ 6\ 5\ 4\ 3\ 1)$,&
  $\pi_{5,5}=(0\ 10)(1\ 3\ 9\ 5)(2\ 4\ 8\ 6)$,&
$\pi_{5,6}=(0\ 4\ 7\ 1\ 3\ 6\ 9\ 10\ 5\ 8)$,\\
$\pi_{6,2}=(0\ 7\ 3\ 10\ 8\ 4)(1\ 6\ 2\ 9\ 5)$,&
$\pi_{11,0}=(1)$.
\end{tabular}
\end{center}

\noindent Let
$E=\{(3,0),(4,0),(4,6),(5,0),(5,1),(5,2),(6,0),(6,1),(6,3),(6,4)\}$
and $M=\{(4,1),(4,2)$, $(4,7),(5,5)\}$. Then for each $(s,t)\in
Adm(10)\setminus(E\cup M)$, $|\pi_{s,t}{\cal B}_1\cap {\cal B}_1|=s$
and $|T(\pi_{s,t}{\cal B}_1\setminus{\cal B}_1)\cap T({\cal
B}_1\setminus\pi_{s,t}{\cal B}_1)|=t$. For each $(s,t)\in M$,
$|\pi_{s,t}{\cal B}_2\cap {\cal B}_1|=s$ and $|T(\pi_{s,t}{\cal
B}_2\setminus{\cal B}_1)\cap T({\cal B}_1\setminus\pi_{s,t}{\cal
B}_2)|=t$.

Now it remains to show that for each $(s,t)\in E$, we have
$(s,t)\not\in Fin(11)$. By Theorem \ref{11non-isomorphic}, there are
exactly $2$ non-isomorphic $(K_4-e)$-designs of order $11$. Thus we
can check all the cases by computer exhaustive search for the fine
triangle intersection numbers of a pair of $(K_4-e)$-designs of
order $11$, i.e., for any permutation $\pi$ on $X$ and for each
$i,j=1,2$, count $|\pi{\cal B}_i\cap {\cal B}_j|$ and $|T(\pi{\cal
B}_i\setminus{\cal B}_j)\cap T({\cal B}_j\setminus\pi{\cal B}_i)|$.
This completes the proof. \qed

\noindent \textbf{Remark}: In this paper, we focus on the
enumerations of $(K_4-e)$-designs of orders $v=6,10,11$. As an
application the fine triangle intersection problem for
$(K_4-e)$-designs of orders $v=6,10,11$ are considered. The
determination of the set $Fin(v)$ is currently being investigated
for any $v\equiv 0,1\ ({\rm mod }\ 5)$ and $v\geq 6$.

\end{document}